\documentclass[12pt]{article}

\usepackage[centertags]{amsmath}
\usepackage{amsfonts}
\usepackage{amssymb}
\usepackage{amsthm}
\usepackage{newlfont}

\newlength{\defbaselineskip}
\newcommand{\setlinespacing}[1]%
           {\setlength{\baselineskip}{#1 \defbaselineskip}}

\theoremstyle{plain}
\newtheorem{thm}{Theorem}[section]

\newtheorem{lem}[thm]{Lemma}

\theoremstyle{definition}
\newtheorem{defn}[thm]{Definition}
\newtheorem{rem}[thm]{Remark}

\numberwithin{equation}{section}

\begin{document}

\newcommand{\N}{\mathbb{N}}
\newcommand{\ol }{\overline}
\newcommand{\ul }{\underline }
\newcommand{\ra }{\rightarrow }
\newcommand{\lra }{\longrightarrow }
\newcommand{\ga }{\gamma }
\newcommand{\st }{\stackrel }
\newcommand{\scr }{\scriptsize }

\title{\Large\textbf{ Polynilpotent Multipliers of Some Nilpotent Products of Cyclic Groups }}
\author{\textbf{Azam Hokmabadi, Behrooz Mashayekhy,}\footnote{This research was in part supported by a grant from IPM (No. 85200037).} \\ \textbf{and Fahimeh Mohammadzadeh} }
\date{ }
\maketitle
\begin{abstract}
In this article, we present an explicit formula for the $c$th
nilpotent multiplier (the Baer invariant with respect to the
variety of nilpotent groups of class at most $c\geq 1$) of the $n$th nilpotent
product of some cyclic groups $G={\mathbb
{Z}}\stackrel{n}{*} \cdots \stackrel{n}{*}{\mathbb
{Z}}\stackrel{n}{*} {\mathbb {Z}}_{r_1}\stackrel{n}{*}
\cdots \stackrel{n}{*}{\mathbb{Z}}_{r_t}$, (m-copies of $\mathbb {Z}$), where $r_{i+1} | r_i$ for
$1 \leq i \leq t-1$ and $c \geq n$ such that $ (p,r_1)=1$ for all
primes $p$ less than or equal to $n$. Also, we compute the
polynilpotent multiplier of the group $G$ with respect to the
polynilpotent variety ${\mathcal N}_{c_1,c_2,\ldots ,c_t}$, where
$c_1 \geq n.$
\end{abstract}
\ \\ \ \\ \ \\ \ \\ \ \\
$2010$ \textit{Mathematics Subject Classification}. 20E10, 20E34, 20F12, 20F18.\\
\textit{Key words and phrases}. Polynilpotent multiplier, nilpotent product, cyclic group.

\newpage
\section{Introduction and Motivation}
\hspace{0.5cm}
Let $G$ be any group with a free presentation $G\cong F/R$, where $F$ is
a free group. Then the Baer invariant of $G$ with respect to the
variety of groups $\mathcal{V}$, denoted by ${\mathcal V}M(G)$,
is defined to be
$${\mathcal V}M(G)=\frac{R \cap V(F)}{[RV^{*}F]},$$
where $V$ is the set of laws of the variety $\mathcal{V}$, $V(F)$
is the verbal subgroup of $F$ and
$$ [R V^{*}F]=\langle v(f_1,\ldots,f_{i-1},f_{i}r,f_{i+1},\ldots,f_k)
 v(f_1,\ldots,f_{i},\ldots,f_k)^{-1} |$$ $$ \hspace{5cm} r\in R, f_i\in F, v
\in V, 1\leq i\leq k, k \in {\mathbb N} \rangle.$$ For example, if
$\mathcal{V}$ is the variety of abelian groups $\mathcal{A}$,
then the Baer invariant of the group $G$ will be $ (R \cap F')/[R,
F],$ which is isomorphic to $M(G)$, the Schur multiplier of $G$ (see [5]).
If $\mathcal{V}$ is the variety of polynilpotent groups of class
row $(c_1,\ldots,c_t)$, $ {\mathcal N}_{c_1,c_2,\ldots,c_t}$, then
the Baer invariant of a group $G$ with respect to this variety, which we call a polynilpotent multiplier, is
as follows:
$${\mathcal N}_{c_1,c_2,\ldots,c_t} M(G)=\frac{R \cap \gamma_{c_t+1}\circ
\ldots \circ \gamma_{c_1+1}(F)}{[R,\ _{c_1}F,\
_{c_2}\gamma_{c_1+1}(F),\ldots,\ _{c_t}\gamma_{c_{t-1}+1}\circ
\ldots \circ \gamma_{c_1+1}(F)]},\ \ \ \ \ \ (1)$$ where
$\gamma_{c_t+1}\circ \cdots \circ
\gamma_{c_1+1}(F)=\gamma_{c_t+1}(\gamma_{c_{t-1}+1}( \cdots
(\gamma_{c_1+1}(F)) \cdots ))$ is the group which is attained from the iterated terms of the lower
central series of $F$. See [4] for the equality $$[R {\mathcal
N^{*}}_{c_1,c_2,\ldots,c_t}F]=[R,\ _{c_1}F,\
_{c_2}\gamma_{c_1+1}(F), \ldots,\ _{c_t}\gamma_{c_{t-1}+1}\circ
\ldots \circ \gamma_{c_1+1}(F)].$$
Note that the Baer invariant of $G$ is always abelian and independent of the
choice of the free presentation for $G$ with respect to a variety ${\mathcal V}$ (see [5]). In particular, if $t=1$
and $c_1=c$, then the Baer invariant of $G$ with respect to the
variety ${\mathcal N}_c$ is called the $c$-nilpotent multiplier
and given by
$${\mathcal N}_cM(G) =\frac{R \cap \gamma
_{c+1}(F)}{[R,\ _cF]}.$$

Determining these Baer invariants of groups is known to be
very useful for classification of groups into isologism classes with respect to the chosen varieties
(see [5]). In 1979, Moghaddam [8] gave a formula for the
$c$-nilpotent multiplier of a direct product of two groups, where
$c+1$ is a prime number or 4. Also, in 1998, Ellis [1] presented
the formula for all $c \geq 1.$ In 1997, Moghaddam and Mashayekhy
[7] presented an explicit formula for the $c$-nilpotent multiplier
of a finite abelian group for every $c\geq1$.\

It is known that the nilpotent product is a generalization of the
direct product. In 1992, Gupta and Moghaddam [2] calculated the
$c$-nilpotent multiplier of the nilpotent dihedral group of class
$n, G_n =\langle x,y | x^2, y^2, [x,y]^{2^{n-1}}\rangle$. It is
routine to verify that $G_n \cong {\mathbb{Z}}_2 \stackrel{n}{*}
{\mathbb{Z}}_2$. In 2003, Moghaddam, Mashayekhy, and Kayvanfar [9]
extended the previous result and calculated the $c$-nilpotent multiplier of $n$th nilpotent products of two
 cyclic groups for $n=$2, 3 and 4 under some conditions. Also, the second author [6]
 gave an implicit formula for the $c$-nilpotent multiplier of a nilpotent product of
 cyclic groups.

In this paper, we first obtain an explicit formula for the
$c$-nilpotent multiplier of the $n$th nilpotent product of some
cyclic groups $G= \underbrace{{\mathbb{Z}}\stackrel{n}{*} \cdots
\stackrel{n}{*}{\mathbb{Z}}}_{m-copies}\stackrel{n}{*}
{\mathbb{Z}}_{r_1}\stackrel{n}{*}\cdots
\stackrel{n}{*}{\mathbb{Z}}_{r_t}$, where $r_{i+1}\mid r_i $ for $1
\leq i \leq t-1$ and $c\geq n$ such that $(p,r_1)=1$ for all
primes $p$ less than or equal to $n$. This result extends the
works of Moghaddam and Mashayekhy [7] and Moghaddam, Mashayekhy
and Kayvanfar [9]. Second, we present an explicit formula for the
polynilpotent multiplier of such a group $G$ with respect to the
polynilpotent variety ${\mathcal N}_{c_1,c_2,\ldots,c_t}$, where
$c_1 \geq n$.

\section{Notation and Preliminaries}

\begin{defn} ([3, \S11.1 and \S12.3]). The basic commutators on the letters $x_{1}, x_{2},\ldots,x_{n},\ldots$ are
defined
as follows: \\
(i) The letters $x_{1}, x_{2},\ldots,x_{n},\ldots $ are basic
commutators of weight one, ordered by setting $x_{i} < x_{j}$, if
$i<j$.\\
(ii) Having defined the basic commutators of weight less than $n$,
the basic commutators of weight $n$ are defined as $c_k=[c_i,c_j]$,
where

(a) \ $c_i$ and $c_j$ are basic commutators and $w(c_i)+w(c_j)=n$,
where $w(c)$ is the weight of $c$ and

(b) \ $c_i>c_j$, and if $c_i=[c_s,c_t]$, then $c_j\geq c_t$.\\
(iii) The basic commutators of weight $n$ follow those of weights
less than $n$. The basic commutators of weight $n$ are ordered
among themselves lexicographically; that is, if $[b_1,a_1]$ and
$[b_2,a_2]$ are basic commutators of weight $n$, then $[b_1,a_1]<
[b_2,a_2]$ if and only if $b_1<b_2$ or $b_1=b_2$
and $a_1<a_2$.
\end{defn}

Basic commutators are special cases of outer commutators. Outer
commutators on the letters $x_{1}, x_{2},\ldots,x_{n},\ldots$ are
defined inductively  as follows:

 The letter $x_i$ is an outer commutator
word of weight one. If $u=u(x_1,\ldots,x_s)$ and
$v=v(x_{s+1},\ldots,x_{s+t})$ are outer commutator words of weights
$s$ and $t$, then
$w(x_1,\ldots,x_{s+t})=[u(x_1,\ldots,x_s),v(x_{s+1},\ldots,x_{s+t})]$
is an
outer commutator word of weight $s+t$ and will be written $w=[u,v]$.
\begin{thm} $([3, \S11.2]).$ Let $F$ be the
free group on $x_{1}, x_{2}, \ldots,x_{d}$, then for all $1 \leq i
\leq n$, $$ \frac{\gamma_{n}(F)}{\gamma_{n+i}(F)} $$ is the free
abelian group, and freely generated by the basic commutators of weights $ n,
n+1,\ldots, n+i-1 $  on $d$ letters.
\end{thm}
\begin{thm} $([3, \S11.2]).$ The number of
basic commutators of weight $n$  on $d$ generators is given by the
following formula: $$ \chi_{n}(d)= \frac{1}{n}\sum_{m|n}\mu(m)
d^{\frac{n}{m}},$$ where $\mu(m)$ is the M\"{o}bius function,
which is defined to be
 $$\mu(m)=\left\{\begin{array}{ll}
 1& if \ m=1,\\
 0& if \ m=p_1^{\alpha_1}\cdots p_k^{\alpha_k} \ \ \  \exists \alpha_{i} > 1,\\
(-1)^s& if \ m=p_1 \cdots p_s, \end{array}\right.$$ where
$p_i$'s are distinct prime numbers.
\end{thm}

Let $G_{i}= \langle x_{i}| x_{i}^{k_{i}} \rangle$, for $i \in I$,
be the cyclic group of order $ k _{i} $ if $ k _{i} > 1 $, and
the infinite cyclic group if $ k _{i} = 0 $. The $n$th nilpotent
product of the family $ \{ G_{i}\} _{i \in I } $ is defined as
follows (see [10]):
$$\prod ^{\stackrel{n}{*}} _{i \in I }G_{i} = \frac{ \prod ^{* } _{i
\in I } G_{i}}{ \gamma _{n+1}(\prod  ^{* } _{i \in I }G_{i})},$$
where $ \prod  ^{* } _{i \in I }G_{i} $ is the free product of the
family $ \{ G_{i} \}_{ i \in I }$.\\
Let $$1 \rightarrow R _{i}= \langle x_{i}^{k_{i}} \rangle
\rightarrow F_{i}= \langle x_{i} \rangle \rightarrow G_{i}
\rightarrow 1 $$ be a free presentation for $ G_{i} $. It is
routine to check that a free presentation for the $n$th nilpotent
product of $ \prod  ^{\stackrel{n}{*}} _{i \in I }G_{i}$ is as
follows (see [9]):
 $$ 1 \rightarrow R=S \gamma _{n+1}(F)\rightarrow F= \prod
^{* } _{i \in I } F_{i} \rightarrow \prod  ^{\stackrel{n}{*}} _{i
\in I }G_{i} \rightarrow 1,$$ where $S=\langle x_{i}^{k_{i}}|i \in
I \rangle ^{F}$. Therefore, if $ c \geq n $, then the $c$-nilpotent
multiplier of $\prod ^{\stackrel{n}{*}} _{i \in I }G_{i}$ is
 $$ {\mathcal N}_{c}M( \prod  ^{\stackrel{n}{*}} _{i
\in I }G_{i})=\frac{R \cap \gamma _{c+1}(F)}{[R,\ _ cF]}=
\frac{\gamma _{c+1}(F)}{[S, \ _ cF]\gamma
_{c+n+1}(F)}=\frac{\gamma _{c+1}(F)}{\rho _{c+1}(S)\gamma
_{c+n+1}(F)},$$ where $\rho _{k}(S)$ is defined inductively by $\rho _{1}(S)=S$ and
 $ \rho _{c+1}(S)=[\rho _{c}(S), F]$.
\begin{lem} If $1\leq i < r$ and $(p,r)=1$
for all primes $p$
less than or equal to $i$, then $r$ divides ${r \choose i}$.
\end{lem}
\begin{proof} Clearly ${r \choose i}=r( \frac{(r-1) \cdots
(r-i+1)}{1 \times 2 \times \cdots \times i})$ is an integer. For
any prime $p\leq i$, $p |(r-1) \cdots
(r-i+1)$, since $p \mid \!\!\!\!/r$. Thus,
$1 \times 2 \times \cdots \times i |(r-1) \cdots (r-i+1)$ and, hence, the result holds.
\end{proof}

The following consequences of the collecting process are
vital in the proof of our main result.
\begin{lem} $([10]).$ Let $x,y$ be any
elements of a given group and let $c_1,c_2,\ldots $ be the sequence of
basic commutators of weights at least two in $x$ and $[x,y]$, in
ascending order. Then
$$ \hspace{3cm} [x^{n},y]=[x,y]^{n}c_1^{f_1(n)}c_2^{f_2(n)} \cdots c_i^{f_i(n)} \cdots,\
\ \hspace{3.3cm} (2)$$where$$ \hspace{3cm} f_i(n)=a_1{n \choose
1}+ a_2{n \choose 2}+ \cdots +a_{w_i}{n \choose w_i},
\hspace{3cm} (3) \
$$with $a_i \in \mathbb{Z}$ and $w_i$ being the weight of $c_i$ in $x$
and $[x,y]$. If the group is nilpotent, then the expression in (2)
gives an identity, and the sequence of commutators terminates.
\end{lem}
\begin{lem} $([10]).$ Let $\alpha$ be a fixed
integer and $G$ a nilpotent group of class at most $n$. If $b_j
\in G$ and $r<n$, then
$$[b_1,..,b_{i-1},b_i^{\alpha},b_{i+1},\ldots ,b_r]=[b_1,\ldots ,b_r]^{\alpha}
c_1^{f_1(\alpha)}c_2^{f_2(\alpha)} \cdots ,$$ where the $c_k$'s are
commutators in $b_1,\ldots ,b_r$ of weight strictly greater than
$r$, and every $b_j$, $1\leq j \leq r$ appears in each commutator
$c_k$, the $c_k$'s listed in ascending order. The $f_i$'s are of the
form (3), with $a_j \in \mathbb{Z}$ and $w_i=(the\ weight\ of\ c_i\ on\
the\ b_i)-(r-1)$.
\end{lem}

\section{Main Results}
\hspace{0.5cm}
Keeping the previous notation, let $k_i=0$, for $1 \leq i \leq m$,
and $k_{m+j}=r_j > 1 $ such that $r_{j+1}|r_j$, for $1 \leq j
\leq t$, then $ \prod  ^{* ^{n}} _{i \in I }G_{i}=
\underbrace{{\mathbb{Z}}\stackrel{n}{*} \cdots
\stackrel{n}{*}{\mathbb{Z}}}_{m-copies}\stackrel{n}{*}
{\mathbb{Z}}_{r_1}\stackrel{n}{*} \cdots
\stackrel{n}{*}{\mathbb{Z}}_{r_t}$. In order to compute the
$c$-nilpotent multiplier of the above group,
we need two technical lemmas.
\begin{lem} With the above notation and
assumption, if $(p, r_1)=1$, for all primes $p$ less than or
equal to $l-i$, then $\rho _{c+i}(S)\gamma _{c+l}(F)/ \rho
_{c+i+1}(S)\gamma _{c+l}(F)$ is the free abelian group with a
basis $D_{i,1}\cup \cdots \cup
D_{i,t}$, where\\
$$D_{i,j}=\{ b^{r_j}\rho_{c+i+1}(S)\gamma _{c+l}(F) |
 \mbox{ b is a basic commutator of weight c+i on} \hspace{0.5cm}$$ $$
 \hspace{3cm} \ x_1,\ldots,x_m,\ldots,x_{m+j}\mbox{ such that } x_{m+j}\mbox{ appears in b}  \},$$
 for $1 \leq i \leq l-1$ and $1 \leq j \leq t$.
\end{lem}
\begin{proof} Using the collecting process (see [3, \S11.1]), one
can easily check that\\ $\rho _{c+i}(S)\gamma _{c+l}(F)/ \rho
_{c+i+1}(S)\gamma _{c+l}(F)$ is generated by all $b'\rho_{c+i+1}(S)\gamma _{c+l}(F)$, where $b'$ belongs to the set of  basic commutators of weight $c+i,$
$\ldots , c+l-1$ on letters $x_1,\cdots , x_m,x_{m+1},\ldots ,
x_{m+t}$ such that one of the $x_{m+1}^{r_1},\ldots , x_{m+t}^{r_t}$
appears in them. It is easy to check that all the above
commutators of weight greater than $c+i$ belong to $\rho
_{c+i+1}(S)$. Now, we show that if $b'$ is one of the
above commutators of weight $c+i$ such that $x_{m+j}^{r_j}$
appears in it, then $$ \hspace{3.5cm} b'\equiv b^{r_j} \pmod{ \rho
_{c+i+1}(S)\gamma _{c+l}(F)}, \hspace{3.5cm} (4) $$ where $b$ is a
basic commutator of weight $c+i$ on $
x_1,\ldots,x_m,\ldots,x_{m+t}$ such that $x_{m+j}$ appears in it.
(Note that $b$ is actually a basic commutator according to the
definition, and $b'$ is the same as $b$, but the letter $x_{m+j}$ with
exponent $r_j$.) In order to prove the above claim, first we use
reverse induction on $k$ $(i+1 \leq k \leq l-1)$ to show that if
$u$ is an outer commutator of weight $c+k$ on
$x_1,\ldots,x_m,\ldots,x_{m+t}$ such that $x_{m+j}$ appears in
$u$, then
$$ \hspace{3.2cm} u^{r_j}\in \rho _{c+i+1}(S) \pmod {\gamma
_{c+l}(F)}.\ \hspace{4cm} (5)$$ Let $k=l-1$ and
$u=[\ldots,x_{m+j},\ldots]$, then clearly $u^{r_j} \equiv
[\ldots,x_{m+j}^{r_j},\ldots]\in \rho _{c+i+1}(S) \pmod{\gamma
_{c+l}(F)}$.

Now, suppose the above property holds for every $k>k'$. We will
show that the property (5) holds for $k'$. Let $u=[u_1,u_2]$ be an
outer commutator of weight $c+k'$ on $x_1,\ldots,x_{m+t}$, where
$x_{m+j}$ appears in $u_1$. Then, by Lemma 2.5, we have
$$u^{r_j} \equiv
[u_1^{r_j},u_2](v_1^{f_1(r_j)} \cdots v_h^{f_h(r_j)})^{-1}
\pmod{\gamma _{c+l}(F)},$$where $v_s$ is a basic commutator of
weight $w_s$ in $u_1$ and $[u_1,u_2]$, $1\leq s\leq h$. Thus, $v_s$ is an outer
commutator of weight greater than $c+k'$ and less than $c+l$ on $
x_1,\ldots,x_m,\ldots,x_{m+t}$ such that $x_{m+j}$ appears in it.
By the hypothesis, since $r_j|r_1$ we have $( p, r_j)=1$ for all
primes $p$ less than or equal to $l-i$. Also, it is easy to see
that $w_s \leq (c+l)-(c+k'-1)=l-k'+1 \leq l-i$. Therefore, by Lemma
2.4, $r_j|f_s(r_j)$, and so, by induction hypothesis, $v_s^{f_s(r_j)}
\in \rho _{c+i+1}(S) \pmod{\gamma _{c+l}(F)}$. Hence, by repeating
the above process, if $u=[\ldots ,x_{m+j},\ldots]$, then $u^{r_j}
\equiv [\ldots,x_{m+j}^{r_j},\ldots]v_1'^{f_1'(r_j)} \cdots v
_h'^{f_h'(r_j)}\in \rho _{c+i+1}(S) \pmod{\gamma _{c+l}(F)}.$
 Now using (5), Lemma 2.6, and some commutator manipulations,
the congruence (4) holds. Therefore, the set $\bigcup
^{t}_{j=1}D_{i,j}$ is a generating set for $\rho _{c+i}(S)\gamma
_{c+l}(F)/ \rho _{c+i+1}(S)\gamma _{c+l}(F)$. On the other hand,
by Theorem 2.2, distinct basic commutators are linearly independent
and, hence, distinct powers of these basic commutators are also linearly independent. Therefore, the set $\bigcup
^{t}_{j=1}D_{i,j}$ is a basis.
\end{proof}
\begin{lem} With the notation and
assumption of the previous lemma, if $(p,r_1)=1$ for all primes $p$
less than or equal to $l-1$, then $$\rho _{c+1}(S)\gamma _{c+l}(F)/
\gamma _{c+l}(F)$$ is the free abelian group with a basis $ \bigcup
_{i=1}^{l-1}( \bigcup ^{t}_{j=1}D_{i,j}).$
\end{lem}
\begin{proof} Put
$$A_i= \frac{\rho _{c+i}(S)\gamma _{c+l}(F)}{
\rho _{c+i+1}(S)\gamma _{c+l}(F)}, B_i=\frac{ \rho _{c+1}(S)\gamma
_{c+l}(F)}{ \rho _{c+i+1}(S)\gamma _{c+l}(F)}.$$ Then, clearly the
following exact sequence exists for $1\leq i \leq l-1$
$$0 \rightarrow A_i \rightarrow B_i \rightarrow B_{i-1}
\rightarrow 0. $$ By Lemma 3.1, $B_1$ is a free abelian group, so the following
exact sequence:
$$0 \rightarrow A_2 \rightarrow B_2 \rightarrow B_1
\rightarrow 0 $$ splits and, hence, $B_2 \cong
A_2 \oplus B_1 $. Also, by Lemma 3.1 every $A_i$ is free abelian,
so by induction, every $B_i$ is free abelian and $$\frac{\rho
_{c+1}(S)\gamma _{c+l}(F)}{\gamma _{c+l}(F)}=B_{l-1} \cong A_{l-1}
\oplus A_{l-2} \oplus \cdots \oplus A_2 \oplus A_1. $$ Now, using
the basis for $A_i$ presented in Lemma 3.1, the result holds.
\end{proof}

Now, we are in a position to state and prove the first main result
of the paper.
\begin{thm} Let $ G=\underbrace{{\mathbb
{Z}}\stackrel{n}{*} \cdots \stackrel{n}{*}{\mathbb
{Z}}}_{m-copies}\stackrel{n}{*} {\mathbb {Z}}_{r_1}\stackrel{n}{*}
\cdots \stackrel{n}{*}{\mathbb{Z}}_{r_t}$ be the $n$th nilpotent
product of some cyclic groups, where $r_{i+1}$ divides $r_i$ for $
1 \leq i \leq t$. If $c \geq n$ and $(p,r_1)=1$ for all primes $p$
less than or equal to $n$, then the $c$-nilpotent multiplier of
$G$ is isomorphic to
$${\mathbb{Z}}^{(d_m)} \oplus {\mathbb{Z}}_{r_1}^{(d_{m+1}-d_{m})}\oplus \cdots \oplus
{\mathbb{Z}}_{r_t}^{(d_{m+t}-d_{m+t-1})},$$ where $d_m=\sum
^{n}_{i=1}\chi _{c+i}(m)$ and ${\mathbb{Z}}_{r_i}^{(d)}$ denotes the
direct sum of $d$ copies of the cyclic group ${\mathbb Z}_{r_i}$.
\end{thm}
\begin{proof} Using the previous notation and assumption, we
have $${\mathcal N}_{c}M( G)= \frac{\gamma _{c+1}(F)}{ \rho
_{c+1}(S)\gamma _{c+n+1}(F)}\cong \frac{\gamma _{c+1}(F)/\gamma
_{c+n+1}(F)}{\rho _{c+1}(S)\gamma _{c+n+l}(F)/ \gamma
_{c+n+1}(F)}.$$
Also, by Theorem 2.2, $\gamma _{c+1}(F)/\gamma
_{c+n+1}(F)$ is a free abelian group with the basis consisting of all basic
commutators of weight $c+1,\ldots,c+n$ on the letters
$x_1,\ldots,x_{m+t}$. 

Now, by considering the basis presented for
$\rho _{c+1}(S)\gamma _{c+n+l}(F)/ \gamma _{c+n+1}(F)$ in Lemma
3.2 and note the points that $D_{i,j}$'s are mutually disjoint and the
number of elements of $D_{i,j}$ is equal to $\chi _{c+i}(m+j)-
\chi _{c+i}(m+j-1)$, the result holds.
\end{proof}

Now the second main result of the paper, which is in turn an
extension of the first one, is as follows:
\begin{thm} Let $G=\underbrace{{\mathbb
{Z}}\stackrel{n}{*} \cdots \stackrel{n}{*}{\mathbb
{Z}}}_{m-copies}\stackrel{n}{*} {\mathbb {Z}}_{r_1}\stackrel{n}{*}
\cdots \stackrel{n}{*}{\mathbb{Z}}_{r_t}$ be the $n$th nilpotent
product of some cyclic groups, where $r_{i+1}$ divides $r_i$, for
$ 1 \leq i \leq t$. If $(p,r_1)=1$ for all primes $p$ less than or
equal to $n$, then the polynilpotent multiplier with class row
$c_1,c_2,\ldots,c_s$ of $G$ is
$${\mathcal N}_{c_1,c_2,\ldots,c_s}M(G)= {\mathbb{Z}}^{(d_m)} \oplus {\mathbb{Z}}_{r_1}^{(d_{m+1}-d_m)}
\oplus \cdots \oplus{\mathbb{ Z}}_{r_t}^{ (d_{m+t}-d_{m+t-1})},$$
where $d_i=\chi_{c_s+1}(\cdots(\chi_{c_2+1}(\sum _{i=1}^{n}\chi
_{c_{1}+i}(m)))\cdots),$ for $c_1 \geq n$
and $c_2, \ldots, c_s \geq 1$ and $ 1 \leq i \leq t$.
\end{thm}
\begin{proof} Let $G$ be a nilpotent group of class $n \leq c_1$
with a free presentation $G=F/R$. Since $\gamma_{c_1+1}(F)\leq
\gamma_{n+1}(F) \leq R$, it gives ${\mathcal N}_{c_1}M(G)=\gamma
_{c_1+1}(F)/[R,\ _{c_1}F]$.\\
Now, we can consider $\gamma
_{c_1+1}(F)/[R,\ _{c_1}F]$ as a free presentation for ${\mathcal
N}_{c_1}M(G)$ and, hence,
$${\mathcal N}_{c_2}M({\mathcal
N}_{c_1}M(G))= \frac{\gamma _{c_2+1}(\gamma _{c_1+1}(F))}{[R,\
_{c_1}F,\ _{c_2}\gamma _{c_1+1}F]}.$$ Therefore, by (1) we have
$${\mathcal N}_{c_1,c_2}M(G)={\mathcal N}_{c_2}M({\mathcal
N}_{c_1}M(G)).$$ By continuing the above process, we can show
that$${\mathcal N}_{c_1,c_2, \ldots,c_t}M(G)={\mathcal
N}_{c_t}M(\cdots {\mathcal N}_{c_2}M({\mathcal
N}_{c_1}M(G))\cdots).$$Using Theorem 3.3, ${\mathcal N}_{c_1}M(G)$
is a finitely generated abelian group of the following form:
 $${\mathbb{Z}}^{(\sum
^{n}_{i=1}\chi _{c_{1}+i}(m))} \oplus {\mathbb{Z}}_{r_1}^{(\sum
^{n}_{i=1}(\chi _{c_{1}+i}(m+1)-\chi _{c_1+i}(m))}\oplus$$ $$ \cdots
\oplus {\mathbb{Z}}_{r_t}^{(\sum ^{n}_{i=1}(\chi _{c_1+i}(m+t)-\chi
_{c_1+i}(m+t-1))}.$$ Now applying Theorem 3.3 with $n=1$, the
result holds.
\end{proof}
\begin{rem} Let $G=
\underbrace{{\mathbb{Z}}\stackrel{n}{*} \cdots
\stackrel{n}{*}{\mathbb{Z}}}_{m-copies}\stackrel{n}{*}
{\mathbb{Z}}_{s_1}\stackrel{n}{*} \cdots
\stackrel{n}{*}{\mathbb{Z}}_{s_t}$
 be the $n$th nilpotent product of some cyclic groups, where the $s_{i}$ are
arbitrary natural numbers, for $ 1 \leq i \leq t$. If $c \geq n$
and $(p,s_i)=1$ for all primes $p$ less than or equal to $n$ and $
1 \leq i \leq t,$ then by a similar proof to Lemmas 3.1 and 3.2 and
Theorem 3.3, one can compute the $c$-nilpotent multiplier of $G$, but
the formula is certainly more complicated than the one in Theorem 3.3.
 For example, if $G={\mathbb{Z}}_{s_1} \stackrel{n}{*} {\mathbb{Z}}_{s_2}
\stackrel{n}{*} {\mathbb{Z}}_{s_3}$, then ${\mathcal N}_{c}M(G)$ is as follows:
$${\mathbb{Z}}_{\alpha}
^{(\sum _{i=1}^{n}\chi _{c+i}(2))}\oplus{\mathbb{Z}}_{\beta}
^{(\sum _{i=1}^{n}\chi _{c+i}(2))}\oplus{\mathbb{Z}}_{\gamma}
^{(\sum _{i=1}^{n}\chi _{c+i}(2))}\oplus {\mathbb{Z}}_{\delta}
^{(\sum _{i=1}^{n}\chi _{c+i}(3)-3 \sum _{i=1}^{n}\chi
_{c+i}(2))},$$ where $\alpha = (s_1,s_2)$, $\beta = (s_2,s_3)$,
$\gamma = (s_1,s_3)$,
$\delta = (s_1,s_2,s_3)$.

Moreover, using the proof of Theorem 3.4 and the above formula twice, we can compute the polynilpotent multiplier with class row
$c_1,c_2$ of $G$ as follows:
$${\mathcal N}_{c_1,c_2}M(G)={\mathbb{Z}}_{\alpha}^{(e_1)}\oplus {\mathbb{Z}}_{\beta}^{(e_1)}\oplus {\mathbb{Z}}_{\gamma}^{(e_1)}\oplus {\mathbb{Z}}_{\delta}^{(e_2)}\oplus {\mathbb{Z}}_{(\alpha , \beta)}^{(e_3)} \oplus {\mathbb{Z}}_{(\alpha , \gamma)}^{(e_3)}\oplus
{\mathbb{Z}}_{(\beta , \gamma)}^{(e_3)}$$ $$\oplus {\mathbb{Z}}_{(\alpha , \delta)}^{(e_4)}\oplus {\mathbb{Z}}_{(\beta , \delta)}^{(e_4)}
\oplus {\mathbb{Z}}_{(\gamma, \delta)}^{(e_4)}\oplus {\mathbb{Z}}_{(\alpha ,\beta , \gamma)}^{(e_5)}\oplus {\mathbb{Z}}_{(\alpha ,\beta , \delta)}^{(e_6)}\oplus {\mathbb{Z}}_{(\beta , \gamma , \delta)}^{(e_6)},$$
where $$e_1= \chi_{c_2+1}(\sum _{i=1}^{n}\chi _{c_1+i}(2)),\ 
e_2=\chi_{c_2+1}(\sum _{i=1}^{n}\chi _{c+i}(3)-3 \sum _{i=1}^{n}\chi_{c+i}(2)),$$
$$ e_3=\chi_{c_2+1}(2\sum _{i=1}^{n}\chi _{c_1+i}(2))-2e_1,\ 
e_4= \chi_{c_2+1}(\sum _{i=1}^{n}\chi _{c+i}(3)-2 \sum _{i=1}^{n}\chi_{c+i}(2))-e_1-e_2,$$
$$ e_5=\chi_{c_2+1}(3\sum _{i=1}^{n}\chi _{c_1+i}(2))-3\chi_{c_2+1}(2\sum _{i=1}^{n}\chi _{c_1+i}(2)),$$
$$ e_6= \chi_{c_2+1}(\sum _{i=1}^{n}\chi _{c+i}(3)- \sum _{i=1}^{n}\chi_{c+i}(2))-\chi_{c_2+1}(2\sum _{i=1}^{n}\chi _{c_1+i}(2))- $$ $$ \chi_{c_2+1}(\sum _{i=1}^{n}\chi _{c+i}(3)-2 \sum _{i=1}^{n}\chi_{c+i}(2)).$$
It seems that the general formula in this case is more complicated than to write!
\end{rem}

Department of Pure Mathematics, P. O. Box 1159-91775, Ferdowsi University of Mashhad, Mashhad, Iran.\\
{\it E-mail address}: hokmabadi$_{-}$ah@yahoo.com\\
\ \ \\
Department of Pure Mathematics, Center of Excellence in Analysis on Algebraic Structures, Ferdowsi University of Mashhad,
P. O. Box 1159-91775, Mashhad, Iran, and Institute for Studies in Theoretical Physics and Mathematics (IPM), Tehran, Iran.\\
{\it E-mail address}: bmashf@um.ac.ir\\
\ \ \\
Department of Pure Mathematics, P. O. Box 1159-91775, Ferdowsi University of Mashhad, Mashhad, Iran.\\
{\it E-mail address}: fa36407@yahoo.com

\begin{thebibliography}{11}
\bibitem{1} G. Ellis, {\it On groups with a finite nilpotent upper
central quotient}, Arch. Math. 70 (1998), 89--96.
\bibitem{2} N. D. Gupta and M. R. R. Moghaddam, {\it Higher Schur
multiplicators of nilpotent dihedral groups}, C. R. Math. Rep.
Acad. Sci. Canada XIV 5 (1992), 225--230.
\bibitem{3} M. Hall, The theory of groups. Macmillan Company, New York, 1959.
\bibitem{4} N. S. Hekster, {\it Varieties of groups and isologisms},
J. Austral. Math. Soc. (Ser. A) 46 (1989), 22--60.
\bibitem{5} C. R. Leedham-Green and S. Mckay, {\it Baer-invariants,
isologism, varietal laws and homology}, Acta Math. 137 (1976), 99--150.
\bibitem{6} B. Mashayekhy, {\it Some notes on the Baer-invariant of a
nilpotent product of groups}, J. Algebra 235 (2001), 15--26.
\bibitem{7} M. R. R. Mogaddam and B. Mashayekhy, {\it Higher Schur
multiplicator of a finite abelian group}, Algebra Colloq.
4 (3) (1997), 317--322.
\bibitem{8} M. R. R. Moghaddam, {\it The Baer-invariant of a direct
product}, Arch. Math. 33 (1979), 504--511.
\bibitem{9} M. R. R. Moghaddam, B. Mashayekhy and S. Kayvanfar,
{\it The higher Schur multiplicator of certain class of groups},
Southeast Asian Bull. Math. 27 (2003), 121--128.
\bibitem{10} R. R. Struik, {\it On nilpotent products of cyclic groups},
Canad. J. Math. 12 (1960), 447--462.
\end{thebibliography}
\end{document}